\documentclass[12pt,a4,oneside]{article}
 \usepackage{afterpage}
 \usepackage{fancyhdr}
  \usepackage{euscript}
 \usepackage[colorlinks]{hyperref}
 \usepackage{lipsum}
 \newcommand\shorttitle{Conformal pointwise-slant Riemannian Maps}
 \newcommand\authors{A. Zaidi, G. Shanker and  J. Yadav}
 
 \usepackage[a4paper, total={6.5in, 9.5in}]{geometry}
 \usepackage{amsmath}
 \usepackage{mathtools}
 \usepackage{amsfonts}
 \usepackage{amssymb}
 \usepackage{amsthm}
 \usepackage{caption}
 \usepackage{subcaption}
 \usepackage{authblk}
 \usepackage[sort]{cite}
 \newtheorem{theorem}{Theorem}[section]

 \newtheorem{example}[theorem]{Example}

 \newtheorem{remark}{\sc Remark}
 \newtheorem{lemma}{\sc Lemma}[section]
 \newtheorem{corollary}{\sc Corollary}[section]
 \newtheorem{definition}{\sc Definition}[section]

 \newcommand{\be}{\begin{eqnarray}}
 	\newcommand{\ee}{\end{eqnarray}}
 \newcommand{\Be}{\begin{eqnarray*}}
 	\newcommand{\Ee}{\end{eqnarray*}}
 \newcommand{\bee}{\begin{equation}}
 	\newcommand{\eee}{\end{equation}}
 \newcommand{\ba}{\begin{array}}
 	\newcommand{\ea}{\end{array}}
 \newcommand{\bl}{\begin{lemma}}
 	\newcommand{\el}{\end{lemma}}
 \newcommand{\bd}{\begin{definition}}
 	\newcommand{\ed}{\end{definition}}
 \newcommand{\bt}{\begin{theorem}}
 	\newcommand{\et}{\end{theorem}}
 \newcommand{\bp}{\begin{proof}}
 	\newcommand{\ep}{\end{proof}}
 \newcommand{\bi}{\begin{itemize}}
 	\newcommand{\ei}{\end{itemize}}
 \newcommand{\br}{\begin{remark}}
 	\newcommand{\er}{\end{remark}}
 \newcommand{\bc}{\begin{corollary}}
 	\newcommand{\ec}{\end{corollary}}
 \newcommand{\bex}{\begin{example}}
 	\newcommand{\eex}{\end{example}}

 \usepackage{chngcntr}
 
 \counterwithin*{equation}{section}
 \counterwithin*{equation}{section}
 
\begin{document}
 	\afterpage{\cfoot{\thepage}}
 	\clearpage
 	\date{}

 	\title{\textbf {Conformal pointwise slant Riemannian maps from or to K\"{a}hler manifolds}}

 	\author{Adeeba Zaidi}
 		\author{Gauree Shanker\thanks{corresponding author, Email: gauree.shanker@cup.edu.in}}
 	 	\author{Jyoti Yadav}
 	\affil{\footnotesize Department of Mathematics and Statistics,
 		Central University of Punjab, Bathinda, Punjab-151 401, India.\\
 		Email:adeebazaidi.az25@gmail.com, gauree.shanker@cup.edu.in*, and sultaniya1402@gmail.com}
 	\maketitle
 		\begin{abstract}
 	  		In this article, we study Conformal pointwise-slant Riemannian maps (\textit{CPSRM}) from or to K\"{a}hler manifolds to or from Riemannian manifolds. To check the existence of such maps, we provide some non-trivial examples. We  derive some important results for these maps. We discuss the integrability and totally geodesicness of the distributions. Further, we investigate the conditions for homotheticity and harmonicity of these maps. Finally, we study some inequalities for these maps.  
 	  	\end{abstract}
 	  	\vspace{0.2 cm}
 	  	
 	  	 	\begin{small}
 	  	 	\end{small}
 	  	 
 	  	 	
 	  	\maketitle
 	  	
 	  	\section{Introduction}
 	  	In differential geometry smooth maps play an important role in study the geometrical properties of a manifold by comparing it with another manifold. Riemannian maps are  the most important type of maps in Riemannian geometry which are the generalization of isometric immersion, Riemannian submersion, and an isometry. In 1992, the notion of Riemannian map was first introduced by Fischer \cite{F}. According to him,  if	$F: (M^m,g_M)\rightarrow (N^n,g_N)$	is a smooth map between smooth finite dimensional Riemannian manifolds $(M,g_M)$ and $(N,g_N)$ such that $0<rankF<min \{m,n\}$ and $F_{*p}: T_pM\rightarrow T_{F(p)}N$ denotes the differential map at $p \in  M$, where $ F(p)  \in  N,$ then $T_pM$ and $T_{F(p)}N$ split
 	  	  	 		 	orthogonally with respect to $ g_M(p)$ and $g_N(F(p))$, respectively, as 
 	  	  	 		 	 	\begin{equation*} 
 	  	  	 		 	 		\begin{split}
 	  	  	 		 	 			T_pM &= kerF_{*p} \oplus (ker F_{*p})^\perp,  \\
 	  	  	 		 	 			&=\mathcal{V}_p\oplus \mathcal{H}_p,
 	  	  	 		 	 		\end{split}	
 	  	  	 		 	 	\end{equation*}
 	  	   		 	 		\begin{equation*}
 	  	   		 	 		T_{F(p)}N = rangeF_{*p}  \oplus  (range F_{*p})^\perp,
 	  	   		 	 	\end{equation*}
 	  	  	 		 	 where $\mathcal{V}_p= kerF_{*p}$ and $\mathcal{H}_p=(ker F_{*p})^\perp$ are vertical and horizontal parts of $T_{p}M$ respectively. The map $F$ is called a Riemannian map at
 	  	  	 		 	 	$p  \in  M,$ if the horizontal restriction
 	  	  	 		 	 	\begin{equation*}
 	  	  	 		 	 		(F_{*p})^h = F_{*p}\ |\ _{\mathcal{H}_p} :\mathcal{H}_p\rightarrow rangeF_{*p}
 	  	  	 		 	 	\end{equation*}
 	  	  	 		 	  is a linear isometry between $(kerF_{*p},g_M\ |_ {kerF_{*p}})$ and $(rangeF_{*p},g_N(y)|_{(rangeF_{*p})}),$ where $y=F(p)$. In other words, $(F_{*p})^h$ satisfies the equation
 	  	  	 		 	 	\begin{equation}\label{fa}
 	  	  	 		 	 		g_N(F_{*} X,F_{*} Y) \ =\ g_M(X,Y),
 	  	  	 		 	 	\end{equation}
 	  	  	 		 	 	for all vector fields $X, Y$ tangent to $\Gamma (kerF_{*p})^\perp $. It can be seen that isometric immersions and Riemannian submersions are particular cases of Riemannian maps with $kerF_{*}=\{0\}$ and $(rangeF_{*})^\perp =\{0\}$ respectively. In 2010, \c{S}ahin \cite{S1} introduced Riemannian maps between almost Hermitian manifolds and Riemannian manifolds. In recent past, many authors have broadly studied various types of Riemannian maps \cite{AK1,S7, G, Z2}. \\    
 	  Moreover,	a smooth map $F : (M^m, g_M) \rightarrow (N^n, g_N)$ between Riemannian manifolds $M$ and $N$ is called a conformal Riemannian map at a point $p\in M$, if there exists a positive function $\lambda(p)$ such that \cite{S7}
 	  	 	  	  	 \begin{equation}
 	  	 	  	  	   g_N(F _*X,F _*Y) = \lambda^2(p)g_M(X,Y) ~~
 	  	 	  	  	            \label{0}
 	  	 	  	  	  \end{equation}
 	  	 	  	  	  for $X,Y \in\varGamma((kerF_{*p})^\perp )$. The function $\lambda(p)$ is called dilation and $\lambda^2(p)$ is the square dilation of $F$ at $p$. $F$ is said to be a conformal Riemannian map, if $F$ is conformal Riemannian  at each point $p\in M.$ It can be seen that for $\lambda=1$, every conformal Riemannian map is a Riemannian map. Further, a conformal Riemannian map $F$ is said to be horizontally homothetic, if the gradient of its dilation $\lambda$ is vertical, i.e., $\mathcal{H}(grad\lambda)=0$ at each point. Conformal Riemannian maps have many applications in various field of science. Therefore, it is very tempting for researchers to investigate different types of conformal Riemannian maps on various structures in complex as well as contact geometry \cite{SY, C2, S4,S3}. Recently, Zaidi et al \cite{Z1} have studied conformal anti-invariant Riemannian maps from or to Sasakian manifolds. \\
 	  	 	  	  	  In this paper, we investigate conformal pointwise-slant Riemannian maps from or to K\"{a}hler manifolds. The paper is divided into four sections. In section $2$, we recall all the basic definitions and terminologies which are needed throughout the paper. In section $3$, we study conformal pointwise-slant Riemannian maps from  K\"{a}hler manifolds to Riemannian manifolds. To show the existence of such maps, we construct an example. We  investigate the integrability of distributions and derive the conditions for horizontal and vertical distributions to be totally geodesic. We establish some results on the homotheticity of the map $F$, we also check the harmonicity of these maps. In section $4$, we investigate  conformal pointwise-slant Riemannian maps from Riemannian manifolds to  K\"{a}hler manifolds and construct an example. We study the integrability of distributions and derive the conditions for horizontal and vertical distributions to be totally geodesic. We drive the condition for the homotheticity and harmonicity of these maps and finally we establish some inequalities for these maps.
 	  	\section{Preliminaries}
 	  	Let $M$ be an even-dimensional manifold. Then a differentiable manifold $M$ is said to be an almost complex manifold, if there exists a linear map $J:TM \rightarrow TM$ satisfying $J^2=-I$ and $J$ is called an almost complex structure of $M$. The tensor field $\mathcal{N}$ of type (1,2), defined by
 	  	\begin{equation}
 	  		\mathcal{N}_J(X,Y)=[JX,JY]-[X,Y]-J[X,JY]-J[JX,Y] \label{ja}
 	  	\end{equation}
 	  	for any $X,Y\in\varGamma (TM)$, is called Nijenhuis tensor field of $J$. If $\mathcal{N}$ vanishes on an almost complex manifold $M$, then $J$ defines a complex structure on $M$ and $M$ is called a complex manifold.  Almost complex manifolds are necessarily orientable. A Riemannian metric $g_M$ on an almost complex manifold $(M,J)$ satisfying 
 	  		\begin{equation}
 	  			g_M(JX,JY)=g_M(X,Y) \label{j}
 	  		\end{equation}
 	  	for all $X,Y \in\varGamma(TM)$, is called an almost Hermitian metric, and the manifold $M$ with Hermitian metric $g_M$ is called an almost Hermitian manifold. If $(\nabla_XJ)Y =0$, for all $X,Y \in\varGamma(TM),$ then $M$ is called a K\"{a}hler manifold \cite{KY}.\\
 	  	Moreover, if $M$ is a K\"{a}hlar manifold, then Riemannian curvature tensor of a complex space form $K(v)$ of constant holomorphic sectional curvature $v$ satisfies \cite{KY}
 	  	 	\begin{equation}\label{space curt}
 	  	 	 \begin{split}
 	  		R_M(Y_1, Y_2, Y_3, Y_4)& = \frac{v}{4}\{g_M(Y_1, Y_4)g_M(Y_2, Y_3)-g_M(Y_1, Y_3)g_M(Y_2, Y_4)\\&+g_M(Y_1, JY_3)g_M(JY_2, Y_4)-g_M(Y_2, JY_3)g_M(JY_1, Y_4)\\&+2g_M(Y_1, JY_2)g_M(JY_3, Y_4)\}
 	  	 	\end{split}
 	  	  \end{equation}
 	 	for vector fields $Y_1, Y_2, Y_3, Y_4\in\Gamma(TK).$\\
 	  	 Further, let $F: (M^m,g_M)\rightarrow (N^n,g_N)$ be a smooth map between smooth finite dimensional Riemannian manifolds, then the differential map $F_{*}$ of $F$ can be viewed as a section of the bundle $Hom(TM, F^{-1}TN)\rightarrow M$, where $F^{-1}TN$ is the pullback bundle whose fibres at $p \in M$ is $(F^{-1}TN)_{p} = T_{F(p)}N$. If the bundle $Hom(TM, F^{-1}TN)$ has a connection $\nabla$ induced from the Levi-Civita connection $\nabla^{M}$ and the pullback connection $\overset{N}{\nabla^{F}}$, then the second fundamental form of $F$ is given by \cite{S7}
 	  	 \begin{equation}  \label{11}
 	  	 (\nabla F_{*})(X,Y) = \overset{N}{\nabla^{F}_{X}}F_{*}Y-F_{*}(\nabla_X^MY) 
 	  	 \end{equation}
 	  	for all $X,Y ~\in \Gamma (TM)$ and $\overset{N}{\nabla^{F}_{X}}F_{*}Y \circ F = \nabla ^{N}_{F_{*}X}F_{*}Y$. \\
 	  Let $F$ be a Riemannian map from a Riemannian manifold $M$ to a Riemannian
 	  manifold $N$. Then we define $\mathcal{T}$ and $\mathcal{A}$ as
 	  \begin{equation}\label{0.1}
 	  \begin{split}
 	  \mathcal{A}_DE &= \mathcal{H}\nabla^
 	  M_{\mathcal{H}D}{\mathcal{V}E} +\mathcal{V}\nabla^M_{\mathcal{H}D}\mathcal{H}E,\\
 	  \mathcal{T}_DE &= \mathcal{H}\nabla^M_{\mathcal{V}D}\mathcal{V}E +\mathcal{V}\nabla^M_{\mathcal{V}D}\mathcal{H}E
 	  \end{split}
 	  \end{equation}
 	  for vector fields $D, E$ on $M$, where $\nabla^
 	  M$ is the Levi-Civita connection of $g_M$. It is also easy to verify that $\mathcal{T}$ is vertical, $\mathcal{T}_D = \mathcal{T}_{\mathcal{V}D}$, and $A$ is horizontal, $\mathcal{A}_D = \mathcal{A}_{\mathcal{H}D}$. On the other hand, from \eqref{0.1} we have \cite{S7}	
 	  	\begin{align}
 	  	 	\nabla _VW \ & = \ \mathcal{T}_VW + \hat{\nabla} _VW, \label{c} \\
 	  	    \nabla _VX \ &= \mathcal{H}\nabla _VX + \mathcal{T} _VX,  \label{d}\\
 	  	 	\nabla _XV \ &= \ \mathcal{A}_XV + \mathcal{V}\nabla _XV, \label{e}\\
 	  	 	\nabla _XY \ &= \ \mathcal{A}_XY + \mathcal{H}\nabla _XY \label{f}
 	  	 \end{align}
 	  	for $X,Y  \in\varGamma((KerF_*)^\perp)$
 	  	 and $V,W \in\varGamma(KerF_*)$, where $\hat{\nabla} _VW =  \mathcal{V}\nabla _VW.$
 	  Also, for any vector field $X$ on $M$ and any section $V$ of $(rangeF{*})^\perp$, we denote by $\nabla ^{F\perp}_{X} V$, the orthogonal projection of $\nabla ^{N}_{X}V $ on $(rangeF_{*})^\perp$, where $\nabla F_{*}^\perp $ is a linear
 	  	connection on $(rangeF_{*})^\perp$ such that $\nabla ^{F\perp} g_{N} = 0.$\\
 	   Further, for a Riemannian map, we have \cite{S7}
 	  	\begin{equation}\label{12}
 	  \nabla ^{N}_{F_{*}X}V=-S_{V}F_{*}X+\nabla^{F\perp}_XV,
 	  	 \end{equation}
 	  	where $S_VF_{*}X$ is the tangential component of  $\nabla ^{N}_{F_{*}X}V$ at $p\in M,~ \nabla ^{N}_{F_{*}X}V(p)\in T_{F(p)}N,~ \\S_{V}F_{*}X(p)\in F_{*p}(T_pM)$ and $\nabla^{F\perp}_XV(p)\in (F_{*p}(T_pM))^\perp.$ It is easy to check that $S_{V}F_{*}X$ is bilinear in $V$ and $F_{*}X$, and $S_{V}F_{*}X$ at $p$ depends only on $V_{p}$ and $F_{*p}X_{p}.$ By direct computations, we can obtain
 	  	\begin{equation}\label{12a}
 	  		g_N(S_VF_*X,F_*Y)=g_N(V,(\nabla F_*)(X,Y))
 	  	 \end{equation}
 	  	for $X,Y\in\Gamma((kerF_*)^\perp)$ and $V\in\Gamma((rangeF_*)^\perp)$.\\
 	  	Moreover, let $F: (M^m,g_M)\rightarrow (N^n,g_N)$ be a conformal submersion \cite{A1}. Then, we have:
 	  	 \begin{equation}\label{Vert. cur}
 	  	 g(R(U, V )W, S) = g(R^{KerF_*}(U, V )W, S) + g(T_UW, T_V S) - g(T_VW, T_US),
 	  	 \end{equation}
 	  	 \begin{equation}\label{hor. cur}
 	  	\begin{split}
 	  	g(R(X, Y)Z, B)& = \frac{1}{\lambda^2}g(R^{(KerF_*^\perp)}(X, Y)Z, B)+\frac{1}{4}\{g(\mathcal{V}[X, Z], \mathcal{V}[Y, B])-g(\mathcal{V}[Y, Z], \mathcal{V}[X, B])\\&+2g(\mathcal{V}[X, Y], \mathcal{V}[Z, B])\} +\frac{\lambda^2}{2}\{g(X, Z)g(\nabla_Ygrad(\frac{
 	  	 		 	 		1}{\lambda^2}), B)\\&- g(Y, Z)g(\nabla_Xgrad(\frac{
 	  	 		 	 		1}{\lambda^2}), B)+g(Y, B)g(\nabla_Xgrad(\frac{
 	  	 		 	 		1}{\lambda^2}), Z)\\&-g(X, B)g(\nabla_Ygrad\frac{1}{\lambda^2}, Z)\}+\frac{\lambda^4}{4}\{(g(X, B)g(Y, Z)\\&-g(Y, B)g(X, Z))||grad(\frac{1}{\lambda^2})||^2\\&+g(X(\frac{1}{\lambda^2})Y-Y(\frac{1}{\lambda^2})X, B(\frac{1}{\lambda^2})Z -Z(\frac{1}{\lambda^2})B)\},		
 	  	 		 	 		\end{split}
 	  	 		 	 	\end{equation}
 	  	 		 	 	where $X, Y,Z,B \in \Gamma(KerF_*)^\perp$ and $U, V,W, S \in \Gamma(KerF_*).$\\
 	  	 		 	  Now, let $F$ be a conformal Riemannian map, then for any $ X,Y \in \Gamma((kerF_{*})^{\perp})$, the second fundamental form $(\nabla F_{*})(X,Y)$ of $F$, is given by \cite{S5} 
 	  	 		\begin{equation}\label{13}
 	  	 	(\nabla F_{*})(X,Y)^{rangeF_{*}} = X(ln\lambda)F_{*}Y + Y(ln\lambda)F_{*}X
 	  	 	- g_{M}(X,Y)F_{*}(grad ln \lambda).
 	  	 		\end{equation}
 	  	 Further, if $(rangeF_{*})^\perp$-component of $(\nabla F_{*})(X,Y) $ is denoted by $(\nabla F_{*})^\perp(X,Y)$, then we can write \cite{G}
 	  	 		\begin{equation} \label{14}
 	  	 		(\nabla F_{*})(X,Y) = (\nabla F_{*})(X,Y)^{rangeF_{*}} + (\nabla F_{*})^\perp(X,Y).
 	  	 		\end{equation}		 	 	
 	  	
 	  	\section{Conformal pointwise slant Riemannian maps \textit{(CPSRM)} from K\"{a}hler manifolds to Riemannian manifolds}
 	  	 In this section, we introduce the notion of conformal pointwise slant Riemannian maps from K\"{a}hler manifolds to Riemannian manifolds, construct an example and discuss the geometry of such maps.
 	  	 
 	  	 \begin{definition}
 	  	Let $F$ be a conformal Riemannian map from a K\"{a}hler manifold $(K,J,g_K)$ to a Riemannian manifold $(L,g_L)$. If for every point $k\in K$, the Wirtinger angle $\theta(X)$ between $JX$ and the space $(kerF_*)_k$ is independent of the choice of $X$, where $X\in \Gamma(kerF_*)$ is a nonzero vector, then $F$ is said to be a conformal pointwise slant Riemannian map. In this case, the angle $\theta$ is regarded as a Function on $K,$ known as slant function of the conformal poinwise slant Riemannian map (\textit{CPSRM}). 
 	  	 \end{definition}
 	  	 Suppose, $F$ be a \textit{CPSRM} from a K\"{a}hler manifold $(K,J,g_K)$ to a Riemannian manifold $(L,g_L)$ with slant function $\theta$, then for any $V\in\Gamma(kerF_*),$ we have
 	  	 \begin{equation}\label{15}
 	  	JV=\phi V+\omega V,
 	  	 \end{equation}
 	  	 where $\phi V\in \Gamma(kerF_*)$ and $\omega V\in \Gamma((kerF_*)^{\perp}).$ Also for any $X\in \Gamma((kerF_*)^\perp),$ we have
 	  	  \begin{equation}\label{16}
 	  	 JX=\mathcal{B}X+\mathcal{C}X,
 	  	 \end{equation} 
 	  	 where $\mathcal{B}X\in\Gamma(kerF_*)$ and $\mathcal{C}X\in \Gamma((kerF_*)^{\perp}).$ Assuming $\mu$ as a orthogonal complementary distribution to $\omega(\Gamma(kerF_*))$ in $\Gamma((kerF_*)^{\perp})$, we can write
 	  	 $$\Gamma((kerF_*)^{\perp})= \omega(\Gamma(kerF_*))\oplus \mu.$$ Further, for a pointwise slant Riemannian map, we have \cite{AK2}
 	  	 \begin{align}
 	  	 (\nabla_V\omega)W&=\mathcal{C}\mathcal{T}_VW-\mathcal{T}_V\phi W, \label{17}\\
 	  	 (\nabla_V\phi)W&=\mathcal{B}\mathcal{T}_VW-\mathcal{T}_V\omega W, \label{18}
 	  	 \end{align}
 	  	 where $\nabla$ is a Levi-Civita coonnection on $K$ and 
 	  	 \begin{align}
 	  	 	  	 (\nabla_V\omega)W&=\mathcal{H}\nabla_V\omega W-\omega \hat{\nabla}_V W, \label{19}\\
 	  	 	  	 (\nabla_V\phi)W&=\hat{\nabla}_V\phi W-\phi \hat{\nabla}_V W \label{20}
 	  	 	  	 \end{align}
 	  	 for $V,W\in \Gamma(kerF_*).$ We say that $\omega$ is parallel with respect to the Levi-Civita
 	  	 connection $\nabla$ on $kerF_*,$ if its covariant derivative with respect to $\nabla$ vanishes,
 	  	 i.e., $(\nabla_V\omega)W = 0$ for $V,W \in \Gamma(ker F_*)$.

\begin{example}
 Consider a Riemannian manifold $(K=\mathbb{R}^4,g_K)$ and a pair of  almost complex structures $\{J_1,J_2\}$ on $K$ satisfying $J_1J_2=-J_2J_1$, where 
 	   	  	 $$J_1(u_1,u_2,u_3,u_4)=(u_3,u_4,-u_1,-u_2)$$
 	   	  	 $$J_2(u_1,u_2,u_3,u_4)=(u_2,-u_1,-u_4,u_3).$$ Let $t:\mathbb{R}^4\rightarrow\mathbb{R}$ be a real-valued function, hence we can define a complex structure 
 	   	  	 $$J_t=(cost)J_1+(sint)J_2$$ on $K$, then $(K,g_K,J_t)$ is an almost complex structure. Again, consider a map $F: (K=\mathbb{R}^4,g_K)\rightarrow (L=\mathbb{R}^4,g_L) $ from a K\"{a}hler manifold $K$ to a Riemannian manifold $L,$ defined by
 	   	  	 $$F(x_1,x_2,x_3,x_4)=(e^{x_1}cosx_3,0,e^{x_1}sinx_3,0),$$ by simple computation we have 
 	   	  	 \begin{equation*}
 	   	  	 kerF_*=span\big\{U=\frac{\partial}{\partial x_2}, V=\frac{\partial}{\partial x_4}\},
 	   	  	 \end{equation*}
 	   	  	\begin{equation*}
 	   	  	(kerF_*)^\perp=span\big\{X=e^{x_1}cosx_3\frac{\partial}{\partial x_1}-e^{x_1}sinx_3\frac{\partial}{\partial x_4},Y=e^{x_1}sinx_3\frac{\partial}{\partial x_1}+e^{x_1}cosx_3\frac{\partial}{\partial x_4}\big \}
 	   	  	\end{equation*} 
 	   	  	and
 	   	   \begin{equation*}
 	   	   	   	  	rangeF_*=span\big\{F_*X=e^{2x_1}\frac{\partial}{\partial y_1},F_*Y=n e^{2x_1}\frac{\partial}{\partial x_3}\big\},
 	   	   	   	  	 \end{equation*}
 	  	hence $F$ is a \textit{CPSRM} from a K\"{a}hler manifold $K$ to a Riemannian manifold $L$ with  $\lambda=e^{x_1}$ and slant function $\theta=t.$	
\end{example} 	  	 
 	  	
 	  	 \begin{lemma}
 	  	 Let $F$ be a Riemannian map from a K\"{a}hler manifold $(K,J,g_K)$ to a Riemannian manifold $(L,g_L)$ with slant function $\theta.$ Then $F$ is a conformal poitwise slant Riemannian map if and only if there exists a constant $\beta\in[-1,0]$ such that $$\phi^2V=\beta V$$ for $V\in \Gamma(kerF_*).$ If $F$ is a conformal slant Riemannian map, then $\beta=-cos^2\theta.$
 	  	 \end{lemma}
 	  	 The proof of the above lemma is exactly same as the proof  for conformal slant Riemannian maps (see \cite{C3}).\\
 	  	 Now, from \eqref{15} and Lemma 3.1, we have the following result.	
 	  	 
 	  	 \begin{lemma}
 	  		Let $F$ be a conformal pointwise slant Riemannian map from a K\"{a}hler manifold $(K,J,g_K)$ to a Riemannian manifold $(L,g_L)$ with slant function $\theta.$ Then, we have
 	  		\begin{align}
 	  		g_K(\phi V,\phi W)&=cos^2\theta g_K( V, W), \label{21}\\g_K(\omega V,\omega W)&=sin^2\theta g_K( V, W) \label{22}
 	  		\end{align}
 	  		for any $V,W\in \Gamma(kerF_*).$
 	  	 \end{lemma} 
 	  	 Also, from \eqref{15} and \eqref{16} we have the following result.
 	  	 \begin{lemma}
 	  	Let $F$ be a \textit{CPSRM} from a K\"{a}hler manifold $(K,J,g_K)$ to a Riemannian manifold $(L,g_L)$. Then, for any $X,Y\in \Gamma((kerF_*)^\perp)$ and $V\in\Gamma(kerF_*)$, we have
 	  	\begin{enumerate}
 	  	\item[$$(i)$$] $g_1(X,\mathcal{C}Y)=-g_1(\mathcal{C}X,Y),$
 	  	\item [$(ii)$] $g_1(\mathcal{C}X,\mathcal{C}Y)=-g_1(X,\mathcal{C}^2Y),$
 	  	\item[$(iii)$] $g_1(X,\mathcal{C}^2Y)=g_1(\mathcal{C}^2X,Y),$
 	  	\item[(iv)] $g_1(X,\omega\phi V)=-g_1(\mathcal{C}X,\omega V).$
 	  	\end{enumerate}
 	  	 \end{lemma} 
 	  	 
 	  	\begin{theorem}
 	  	Let $F$ be a \textit{CPSRM} from a K\"{a}hler manifold $(K,J,g_K)$ to a Riemannian manifold $(L,g_L)$. If $\omega$ is parallel with respect to $\overset{K}{\nabla}$ on $kerF_*,$ then we have $$\mathcal{T}_{\phi V}\phi V=-cos^2\theta \mathcal{T}_{V}V,$$
 	  	where $V\in \Gamma(kerF_*)$ and $\theta$ is a slant function of \textit{CPSRM}.
 	  		\end{theorem}
 	  	\begin{proof}
 	  	Let $V,W\in \Gamma(kerF_*)$ and $\omega$ is parallel with respect to $\overset{K}{\nabla}$ on $kerF_*,$ from \eqref{17}, we have $$\mathcal{C}\mathcal{T}_VW=\mathcal{T}_V\phi W,$$
 	  	interchanging $V$ and $W$ in above equation, subtracting the resultant from above equation, we get
 	  	 \begin{equation}\label{22a}
 	  	 \mathcal{T}_V\phi W=\mathcal{T}_W\phi V.
 	  	 \end{equation}
 	Putting $W=\phi V$ and using lemma 3.1 in \eqref{22a}, we have
 	 \begin{equation}\label{23}
 	 	  	 \mathcal{T}_{\phi V}\phi V=-\mathcal{T}_Vcos^2\theta  V.
 	 	  	 \end{equation} 
 	 	  	 From \eqref{c}, we can write 
 	 	  	\begin{equation}\label{0.2}
 	 	  	\begin{split}
 	 	  	\mathcal{T}_Vcos^2\theta  V&=cos^2\theta\mathcal{H}\nabla_V V-\mathcal{H}(sin2\theta V(\theta)V)\\
 	 	  	 &=cos^2\theta	\mathcal{T}_VV. 
 	 	  	\end{split}	 	  		
 	 	  	\end{equation}
 	 	 Hence, from \eqref{23} and \eqref{0.2}, we get required result. 
 	  	\end{proof}
 	  	
 	  	\begin{theorem}
 	  	Let $F$ be a \textit{CPSRM} from a K\"{a}hler manifold $(K,J,g_K)$ to a Riemannian manifold $(L,g_L)$ with slant function $\theta$. Then, any two of the following assertion imply the third one
 	  	\begin{enumerate}
 	  	\item [$(i)$] $(kerF_*)^\perp$ is integrable,
 	  	\item [$(ii)$] for any $X,Y\in \Gamma((kerF_*)^\perp)$ and $ V\in\Gamma (ker F_*),$ 
 	   \begin{equation*}
 	  	\begin{split}
 	  	g_L\big(\overset{L}{\nabla^{F}_{X}}F_*(\omega\phi V),F_*Y\big)-g_L\big(\overset{L}{\nabla^{F}_{Y}}F_*(\omega\phi V),F_*X\big)&=g_L\big(\overset{L}{\nabla^{F}_{X}}F_*(\omega V),F_*(\mathcal{C}Y)\big)\\&
 	  	-g_L\big(\overset{L}{\nabla^{F}_{Y}}F_*(\omega V),F_*(\mathcal{C}X)\big),
 	  	\end{split}
 	  	\end{equation*}
 	  	\item [$(iii)$] $F$ is a horizontally homothetic map.
 	  	\end{enumerate}
 	  	\end{theorem}	
 	  \begin{proof}
 	  Let $X,Y\in \Gamma((kerF_*)^\perp)$ and $ V\in\Gamma (ker F_*)$, then we have
 	  \begin{equation}\label{0.3}
 	  g_K([X,Y],V)=g_K(\overset{K}{\nabla} _XY-\overset{K}{\nabla}_YX,V).
 	  \end{equation}
 	  Since $K$ is a K\"{a}hler manifold, from \eqref{15} and \eqref{0.3}, we have
 	  \begin{equation*}
 	  \begin{split}
 	  g_K([X,Y],V)&=g_K(\overset{K}{\nabla}_X\phi^2V+\overset{K}{\nabla}_X\omega\phi V,Y)-g_K(\overset{K}{\nabla}_X\omega V,JY)\\&
 	  -g_K(\overset{K}{\nabla}_Y\phi^2V+\overset{K}{\nabla}_Y\omega\phi V,X)+g_K(\overset{K}{\nabla}_Y\omega V,JX),
 	  \end{split}
 	  \end{equation*}
 	  using the property of conformal map and lemma 3.1, above equation can be written as
 	   \begin{equation}\label{24}
 	   	  \begin{split}
 	   	 sin^2\theta g_K([X,Y],V)&=g_K(sin2\theta X(\theta) V,Y)+\frac{1}{\lambda^2}\big(g_L(F_*(\overset{K}{\nabla}_X\omega\phi V),F_*Y)\\&-g_L(F_*(\overset{K}{\nabla}_Y\omega\phi V),F_*X)-g_L(F_*(\overset{K}{\nabla}_X\omega V),F_*(\mathcal{C}Y))
 	   	 \\&+g_L(F_*(\overset{K}{\nabla}_Y\omega V),F_*(\mathcal{C}X))\big).
 	   	  \end{split}
 	   	  \end{equation}
 	   	  Since $F$ is a conformal Riemannian map, using \eqref{11}, \eqref{13} and \eqref{14}, we get
 	   	  \begin{equation*}
 	   	  	   	  \begin{split}
 	   	  	   	 sin^2\theta g_K([X,Y],V)&=\frac{1}{\lambda^2}\Big(g_L\big(-(\nabla F_*)(X,\omega\phi V)-X(ln\lambda)F_*(\omega\phi V)-\omega\phi V(ln\lambda )F_*X\\&+g_K(X,\omega\phi V)F_*(grad(ln\lambda))+\overset{L}{\nabla_X^F}F_*(\omega\phi V),F_*Y\big)\\&-g_L\big(-(\nabla F_*)(Y,\omega\phi V)-Y(ln\lambda)F_*(\omega\phi V)-\omega\phi V(ln\lambda )F_*Y\\&+g_K(Y,\omega\phi V)F_*(grad(ln\lambda))+\overset{L}{\nabla_Y^F}F_*(\omega\phi V),F_*X\big)\\&-g_L\big(-(\nabla F_*)(X,\omega V)-X(ln\lambda)F_*(\omega V)-\omega V(ln\lambda )F_*X\\&+g_K(X,\omega V)F_*(grad(ln\lambda))+\overset{L}{\nabla_X^F}F_*(\omega V),F_*({\mathcal{C}Y})\big)\\&+g_L\big(-(\nabla F_*)(Y,\omega V)-Y(ln\lambda)F_*(\omega V)-\omega V(ln\lambda )F_*Y\\&+g_K(Y,\omega V)F_*(grad(ln\lambda))+\overset{L}{\nabla_Y^F}F_*(\omega V),F_*({\mathcal{C}X})\big)\Big).
 	   	  	   	  \end{split}
 	   	  	   	  \end{equation*}
 	  After simplifying the above equation and using Lemma 3.3, we get
 	  \begin{equation}\label{25}
 	  \begin{split}
 	 sin^2\theta g_K([X,Y],V)&=3X(ln\lambda)g_K(\mathcal{C}Y,\omega V)-3Y(ln\lambda)g_K(\mathcal{C}X,\omega V)\\&+2(\omega V)(ln\lambda)g_K(X,\mathcal{C}Y)+\mathcal{C}X(ln\lambda)g_K(Y,\omega V)-\mathcal{C}Y(ln\lambda)g_K(X,\omega V)\\&
 	 +\frac{1}{\lambda^2}\Big(g_L\big(\overset{L}{\nabla_X^F}F_*(\omega\phi V),F_*Y\big)-g_L\big(\overset{L}{\nabla_Y^F}F_*(\omega\phi V),F_*X\big)\\&-g_L\big(\overset{L}{\nabla_X^F}F_*(\omega V),F_*({\mathcal{C}Y})\big)+g_L\big(\overset{L}{\nabla_Y^F}F_*(\omega V),F_*({\mathcal{C}X})\big)\Big).
 	  \end{split}
 	  \end{equation} 
 	  Now, assuming assertions $(i)$ and $(ii)$ are satisfied by \eqref{25}, and taking $X=Y$, we have
 	  \begin{equation}
 	 g_K(\omega V,\mathcal{H}grad(ln\lambda))g_K(X,\mathcal{C}X)=0,
 	  \end{equation}
 	  which is possible only if $\mathcal{H}grad(ln\lambda)=0,$ this implies $(iii)$. Similarly, one can easily show that assertions $(ii)$ and $(iii)$, imply $(i)$ and assertions $(i)$ and $(iii)$ imply $(ii)$. Hence, the theorem.
 	  \end{proof}
 	  \begin{theorem}
 		Let $F$ be a \textit{CPSRM} from a K\"{a}hler manifold $(K,J,g_K)$ to a Riemannian manifold $(L,g_L)$. Then vertical distribution $kerF_*$ defines a totally geodesic foliation on $K$ if and only if
 		\begin{equation}\label{26}
 		\lambda^2g_K(\mathcal{T}_V\mathcal{B}X,\omega W)=g_L\big((\nabla F_*)(V,\omega\phi W),F_*X\big)-g_L\big((\nabla F_*)(V,\omega W),F_*(\mathcal{C}X)\big),
 		\end{equation}
 	where $V,W\in\Gamma(kerF_*)$ and $X\in\Gamma((kerF_*)^\perp).$	
 	  \end{theorem}
 	  \begin{proof}
 	  Let $V,W\in\Gamma(kerF_*)$ and $X\in\Gamma((kerF_*)^\perp).$ Since $K$ is a K\"{a}hler manifold, from \eqref{15}, we have
 	  \begin{equation}\label{27}
 g_K(\overset{K}{\nabla} _V W,X)=g_K(\overset{K}{\nabla}_V(\phi W+\omega W),JX),
 	  \end{equation}
 	  using \eqref{15}, \eqref{16} and lemma 3.1 in \eqref{27}, we get
 	  \begin{equation*}
 	\begin{split}
 	 g_K(\overset{K}{\nabla} _V W,X)&=-g_K(\overset{K}{\nabla}_Vcos^2\theta W,X)-g_K(\overset{K}{\nabla}_V\omega \phi W,X)+g_K(\overset{K}{\nabla}_V\omega W,JX),\\
 	 sin^2\theta g_K(\overset{K}{\nabla} _V W,X)&=-sin2\theta g_K(V(\theta) W,X)-g_K(\overset{K}{\nabla}_V\omega \phi W,X)+g_K(\overset{K}{\nabla}_V\omega W,\mathcal{B}X+\mathcal{C}X), 
 	\end{split}
 	  \end{equation*}
 	further, using the condition of conformality and \eqref{c} in above equation, we have
 	 \begin{equation}\label{28}
 	 	\begin{split}
 	 	 sin^2\theta g_K(\overset{K}{\nabla} _V W,X)&=-g_K(\mathcal{T}_V\mathcal{B}X,\omega W)+\frac{1}{\lambda^2}\Big(g_L\big(F_*(\overset{K}{\nabla}_V\omega W),F_*(\mathcal{C}X)\big)\\&-g_L\big(F_*(\overset{K}{\nabla}_V\omega \phi W),F_*X\big)\Big). 
 	 	\end{split}
 	 	  \end{equation}
 	 If $kerF_*$ defines a totally geodesic foliation on $K$, then from \eqref{11} and \eqref{28}, we have \eqref{26}. This completes the proof.	    
 	  \end{proof}
 	  \begin{corollary}
 	 	Let $F$ be a \textit{CPSRM} from a K\"{a}hler manifold $(K,J,g_K)$ to a Riemannian manifold $(L,g_L)$. If $kerF_*$ defines totally geodesic foliation on $K$, then
 	 	\begin{equation*}
 	 g_K(\mathcal{T}_V\mathcal{B}X,\omega V)=2V(ln\lambda)g_K(\omega\phi W,X).
 	 	\end{equation*}
 	 	 \end{corollary}
 	 	\begin{proof}
 	  Let $V,W\in\Gamma(kerF_*)$ and $X\in\Gamma((kerF_*)^\perp),$ then from theorem 3.3, we have
 	  	\begin{equation*}
 	   		\lambda^2g_K(\mathcal{T}_V\mathcal{B}X,\omega W)=g_L\big((\nabla F_*)(V,\omega\phi W)^{rangeF_*},F_*X\big)-g_L\big((\nabla F_*)(V,\omega W)^{rangeF_*},F_*(\mathcal{C}X)\big),
 	   		\end{equation*}
 	   	using \eqref{13} in above equation, we get
 	   	\begin{equation}
 	   		\lambda^2g_K(\mathcal{T}_V\mathcal{B}X,\omega W)=V(ln\lambda)\Big(g_L\big(F_*(\omega\phi W),F_*X\big)-g_L\big(F_*(\omega W),F_*(CX)\big)\Big).
 	   \end{equation}
 	  Applying lemma 3.3 in above equation, we get required result. 	
 	 	\end{proof}
 	 
 	  \begin{theorem}
 	 	Let $F$ be a \textit{CPSRM} from a K\"{a}hler manifold $(K,J,g_K)$ to a Riemannian manifold $(L,g_L)$. Then, any two of the following assertions imply the third one
 	 	\begin{enumerate}
 	 	\item [$(i)$] $(kerF_*)^\perp$ defines a totally geodesic foliation on $K$.
 	 	\item[$(ii)$]  $\lambda$ is constant on $(kerF_*)^\perp$.
 	 	\item [$(iii)$]	\begin{equation}\label{27a}
 	 	 	 	\lambda^2g_K(\mathcal{A}_X\mathcal{B}Y,\omega V)=g_L\big(\overset{L}{\nabla_X^F}F_*(\omega V),F_*(\mathcal{C}Y)\big)-g_L\big(\overset{L}{\nabla_X^F}F_*(\omega\phi V),F_*Y\big),
 	 	 	 	\end{equation}
 	 	 	 	where $X,Y\in\Gamma((kerF_*)^\perp$ and $V\in\Gamma(kerF_*).$
 	 	\end{enumerate}  
 	  \end{theorem}
 	  \begin{proof}
 	 Let $X,Y\in\Gamma((kerF_*)^\perp$ and $V\in\Gamma(kerF_*)$ and $K$ is a K\"{a}hler manifold, then from \eqref{15}, we have
 	 \begin{equation}\label{28a}
 g_K(\overset{K}{\nabla}_XY,V)=-g_K(\overset{K}{\nabla}_X\phi V+\overset{K}{\nabla}_X\omega V,JY),
  \end{equation}
 from \eqref{16}, \eqref{28a} and lemma 3.1, we get
  \begin{equation}\label{29}
 sin^2\theta g_K(\overset{K}{\nabla}_XY,V)=g_K(sin2\theta X(\theta)V,Y)+g_K(\overset{K}{\nabla}_X\omega\phi V,Y)-g_K(\overset{K}{\nabla}_X\omega V,\mathcal{B}Y)-g_K(\overset{K}{\nabla}_X\omega V,\mathcal{C}Y).
  \end{equation}
  Since $F$ is a conformal Riemannian map, from \eqref{e} and \eqref{29}, we have
    \begin{equation*}
    sin^2\theta g_K(\overset{K}{\nabla}_XY,V)=g_K(\mathcal{A}_X\mathcal{B}Y,\omega V)+\frac{1}{\lambda ^2}\Big(g_L\big(F_*(\overset{K}{\nabla}_X\omega\phi V),F_*Y\big)-g_L\big(F_*(\overset{K}{\nabla}_X\omega V),F_*(\mathcal{C}Y)\big)\Big),
     \end{equation*}
    using \eqref{11}, \eqref{13} and \eqref{14}, we have
      \begin{equation}\label{30}
        \begin{split}
         sin^2\theta g_K(\overset{K}{\nabla}_XY,V)&=g_K(\mathcal{A}_X\mathcal{B}Y,\omega V)+\frac{1}{\lambda ^2}\Big(g_L\big(\overset{L}{\nabla^F_X}F_*(\omega\phi V)-X(ln\lambda)F_*(\omega\phi V)-(\omega\phi V)(ln\lambda)F_*X\\&+g_K(X,\omega\phi V)F_*(grad(ln\lambda)),F_*Y\big)-g_L\big(\overset{L}{\nabla^F_X}F_*(\omega V)-X(ln\lambda)F_*(\omega V)\\&-\omega V(ln\lambda)F_*X+g_K(X,\omega V)F_*(grad(ln\lambda)),F_*(\mathcal{C}Y)\big)\Big).
        \end{split}
          \end{equation}
          Assuming, assertions $(i)$ and $(ii)$ are true, from \eqref{30}, we have $(iii)$. Similarly, if assertions $(ii)$ and $(iii)$ are true, then from \eqref{30}, we get $(i)$. Further, if assertions $(i)$ and $(iii)$ are true, taking $X=Y$ and using lemma 3.3 in \eqref{30}, we have
          \begin{equation*}
       (\omega\phi V)(ln\lambda)g_K(X,X)+X(ln\lambda)g_K(\omega\phi V,X)+\mathcal{C}X(ln\lambda)g_K(X,\omega V)=0,
       \end{equation*}
       which implies that
        \begin{align*}
             g_K (\omega\phi V,grad(ln\lambda))&=0,\\
             g_K(X,grad(ln\lambda))&=0,\\
             g_K(\mathcal{C}X,grad(ln\lambda))&=0.
              \end{align*}
              This is possible if and only if $\lambda$ is constant on $(kerF_*)^\perp$. Hence the theorem.
 	  \end{proof}
 	  \begin{theorem}
 	  	Let $F$ be a \textit{CPSRM} from a K\"{a}hler manifold $(K,J,g_K)$ to a Riemannian manifold $(L,g_L)$ with $\theta$ as a slant function. Then $F$ is harmonic if and only if $\omega$ is parallel and $\lambda$ is constant on $(kerF_*)^\perp$.
 	  \end{theorem}
 	  \begin{proof}
 	Consider a canonical orthogonal frame $e_1, sec\theta \phi e_1,e_2, sec\theta \phi e_2,...e_r, sec\theta \phi e_r,csc\theta \omega e_1,...,\\
 	csc\theta \omega e_{2r},\tilde{e}_1,...,\tilde{e}_s$ such that $\{e_1, sec\theta \phi e_1,e_2, sec\theta \phi e_2,...e_r, sec\theta \phi e_r\}$ is an orthonormal basis of $kerF_*$ and $\{\tilde{e}_1,...,\tilde{e}_s\}$ is of $\mu$. Then the map $F$ is said to be harmonic if and only if 
 	\begin{equation}\label{31}
\begin{split}
trace|_{kerF_*} \Big \{&\sum_{i=1}^{r}\big((\nabla F_*)(e_i,e_i)+sec_2\theta(\nabla F_*)(\phi e_i,\phi e_i)\big)\Big\}\\&
+ trace|_{(kerF_*)^\perp}\Big\{csc^2\theta \sum_{i=1}^{2r}(\nabla F_*)(\omega e_i,\omega e_i)+ \sum_{j=1}^{s}(\nabla F_*)(\tilde{e}_j,\tilde{e}_j)\Big\}=0.
\end{split}
 	\end{equation}
 	Since $K$ is a K\"{a}hler manifold and $F$ ia a \textit{CPSRM}, from \eqref{c} and \eqref{11}, we have
 	\begin{equation}\label{32}
 	\sum_{i=1}^{r}\big((\nabla F_*)(e_i,e_i)+sec_2\theta(\nabla F_*)(\phi e_i,\phi e_i)\big)=-\sum_{i=1}^{r}F_*(\mathcal{T}_{e_i}e_i+sec^2\theta \mathcal{T}_{\phi e_i}\phi e_i).
 	\end{equation}
 	Further, from \eqref{13}, \eqref{14} and lemma 3.2, we get
 	\begin{equation}\label{33}
\begin{split}
 csc^2\theta \sum_{i=1}^{2r}(\nabla F_*)(\omega e_i,\omega e_i)+ \sum_{j=1}^{s}(\nabla F_*)(\tilde{e}_j,\tilde{e}_j)&=csc^2\theta \sum_{i=1}^{2r}\big((\nabla F_*)^\perp (\omega {e_i},\omega {e_i})\\&
 +2g_K(\omega e_i,grad(ln\lambda))F_*(\omega e_i)\big)\\&
 +\sum_{j=1}^{s}\big((\nabla F_*)^\perp(\tilde{e}_j,\tilde{e}_j)\\&
 +2g_K(\tilde{e}_j,grad(ln\lambda))F_*(\tilde{e}_j)\big)\\&
 -(2r +s)F_*(grad ln\lambda),
\end{split}
 	\end{equation}
 after simplifying \eqref{33}, we get
 	\begin{equation}\label{34}
 \begin{split}
  csc^2\theta \sum_{i=1}^{2r}(\nabla F_*)(\omega e_i,\omega e_i)+ \sum_{j=1}^{s}(\nabla F_*)(\tilde{e}_j,\tilde{e}_j)&=csc^2\theta \sum_{i=1}^{2r}\big((\nabla F_*)^\perp (\omega {e_i},\omega {e_i})\\&
  +\sum_{j=1}^{s}\big((\nabla F_*)^\perp(\tilde{e}_j,\tilde{e}_j)\\&
    +(4-2r-s)F_*(grad ln\lambda),
 \end{split}
  	\end{equation}
 	since, $F$ is harmonic, from \eqref{31}, \eqref{32}, \eqref{34} and theorem 3.2, we obtain the required result.
 	  \end{proof}

 	 \begin{theorem}
 	  	Let $F$ be a CPSRM from a complex space form $(K(v), g_K)$ to a Riemannian manifold $(L, g_L)$ with $(rangeF_*)^\perp=\{0\}$. Then
 	  	\begin{equation}
 	  Ric^{(kerF_*)}(U)\leq \frac{v}{4}(2r-1+3cos^2\theta)g_K(U, U )+2rg(T_U U, H),		
 	  	\end{equation}
 	 where $U\in\Gamma(kerF_*)$, $H$ is mean curvature vector field, $v$ is constant holomorphic sectional curvature and $dim(kerF_*)=2r$. The equality holds if and only if the fibers are totally geodesic.
 	  \end{theorem}
 	  \begin{proof}
 	   
 	   	  Let $F: (K,g_K, J)\rightarrow (L,g_L)$ be a \textit{CPSRM} with $(rangeF_*)^\perp = \{0\}.$  For every point $p \in\Gamma(TK)$, let $E_1, ..., E_{2r},csc\theta \omega e_1,...,
 	   	   	csc\theta \omega e_{2r},\tilde{e}_1,...,\tilde{e}_s$ be an orthonormal basis of $T_pK(v)$ such that $kerF_* = span\{E_1, ..., E_{2r}\},$ and $ (kerF_*)^\perp = span\{csc\theta\omega e_1, ..., csc \theta\omega e_{2r},\tilde{e}_1,...,\tilde{e}_s\},$ then for any $U, V,W, S \in \Gamma(KerF_*),$
 	   	  using \eqref{Vert. cur}, we have
 	   	  \begin{equation}\label{jy}
 	   	  	g(R^{KerF_*}(U, V )W, S) = g(R_K(U, V )W, S) - g(T_UW, T_V S) + g(T_VW, T_US).
 	   	  \end{equation}
 	   	  Further, from \eqref{space curt} and \eqref{jy} , we get 
 	   	  \begin{equation}
 	   	  	\begin{split}
 	   	  	g(R^{KerF_*}(U, V )W, S)&=  \frac{v}{4}\Big\{g_K(U, S)g_K(V, W)-g_K(U, W)g_K(V, S)\\&+g_K(U, JW)g_K(JV, S)-g_K(V, JW)g_K(JU, S)\\&+2g_K(U, JV)g_K(JW, S)\Big\}.		
 	   	  	\end{split}
 	   	  \end{equation}
 	   	  Putting $U=S$ and $V=W= E_i, i=1,...,2r$ in above equation, for any vertical vector $U$, we have
 	   	  \begin{equation}\label{Ricci vert}
 	   	  Ric^{(KerF_*)}(U) = \frac{v}{4}(2r-1+3cos^2\theta)g_K(U, U)+2rg(T_U U, H)-g(T_U E_i, T_{E_i}U).
 	   	  \end{equation}
 	   	  Hence, from above equation we get the required result.
 	  \end{proof}
 	 
 	  \begin{theorem}
 	  		Let $F$ be a CPSRM from a complex space form $(K(v), g_K)$ to a Riemannian manifold $(L, g_L)$ with $(rangeF_*)^\perp=\{0\}$. Then
 	  	\begin{equation}
 	  		\begin{split}
 	    \frac{1}{\lambda^2}g(Ric^{(KerF_*^\perp)}(X)&\leq \dfrac{v}{4}\Big\{(2r+s+2)||X||^2+3g_K(\omega\mathcal{B}X,X)\Big\}-\frac{1}{4}||3\mathcal{V}[X,X_n]||^2\\&-\frac{\lambda^2}{2}\{2g_K(X, X_n)g_K(\nabla_Xgrad\frac{1}{\lambda^2}, X_n)-(2r+s)g_K(\nabla_{X}grad\frac{1}{\lambda^2}, X)\\&
 	     	   	  -||X||^2g_K(\nabla_{X_n}grad\frac{1}{\lambda^2}, X_n)\big\}-\frac{\lambda^4}{4}\Big\{\big((2r+s)||X||^2\\&
 	     	   	  -(g_K(X, X_n))^2\big)||grad\frac{1}{\lambda^2}||^2\Big\},
 	  		\end{split}
 	  	\end{equation}
 	 where $X\in\Gamma(KerF_*)^\perp$, $\{X_n=X_i+X_j\}_{i=1,...,2r,j=1,...,s}$ is orthonormal basis for $(kerF_*)^\perp$,  $v$ is constant holomorphic sectional curvature and $dim(kerF_*)^\perp=2r+s.$ The equality holds if and only if $F$ is conformal homothetic map.
 	  \end{theorem}
 	  \begin{proof}
 	 Let $F: (K,g_K, J)\rightarrow (L,g_L)$ be a \textit{CPSRM} with $(rangeF_*)^\perp = \{0\}$ and $e_1, sec\theta \phi e_1,e_2,\\ sec\theta \phi e_2,...e_r, sec\theta \phi e_r,csc\theta \omega e_1,...,
 	 	csc\theta \omega e_{2r},\tilde{e}_1,...,\tilde{e}_s$ be a canonical orthogonal frame such that $\{e_1, sec\theta \phi e_1,e_2, sec\theta \phi e_2,...e_r, sec\theta \phi e_r\}$ is an orthonormal basis of $kerF_*$ and $\{\tilde{e}_1,...,\tilde{e}_s\}$ is of $\mu$. Let $X, Y,Z,B \in \Gamma(KerF_*)^\perp$ and $v$ be the constant holomorphic sectional curvature, then from \eqref{hor. cur} and \eqref{space curt}, we have 
 	   	  \begin{equation}
 	   	  	\begin{split}
 	   	   \frac{1}{\lambda^2}g(R^{(KerF_*^\perp)}(X, Y)Z, B)& = \frac{v}{4}\Big\{g_K(X, B)g_K(Y, Z)-g_K(X, Z)g_K(Y, B)\\&+g_K(X, JZ)g_K(JY, B)-g_K(Y, JZ)g_K(JX, B)\\&+2g_K(X, JY)g_K(JZ, B)\Big\} 	-\frac{1}{4}\Big\{g(\mathcal{V}[X, Z], \mathcal{V}[Y, B])\\&-g(\mathcal{V}[Y, Z], \mathcal{V}[X, B])+2g(\mathcal{V}[X, Y], \mathcal{V}[Z, B])\Big\}\\& -\frac{\lambda^2}{2}\Big\{g(X, Z)g(\nabla_Ygrad(\frac{
 	   	  	1}{\lambda^2}), B)- g(Y, Z)g(\nabla_Xgrad(\frac{
 	   	  	1}{\lambda^2}), B)\\&+g(Y, B)g(\nabla_Xgrad(\frac{
 	   	  	1}{\lambda^2}), Z)-g(X, B)g(\nabla_Ygrad\frac{1}{\lambda^2},Z)\Big\}\\&-\frac{\lambda^4}{4}\Big\{\big(g(X, B)g(Y, Z)-g(Y, B)g(X, Z)\big)||grad(\frac{1}{\lambda^2})||^2\\&+g\Big(X(\frac{1}{\lambda^2})Y-Y(\frac{1}{\lambda^2})X, B(\frac{1}{\lambda^2})Z -Z(\frac{1}{\lambda^2})B\Big)\Big\}.		
 	   	  	\end{split}
 	   	  \end{equation}
 	   	  By using, $X=B$ and $Y=Z=csc\theta\omega e_i+\tilde{e}_j=X_i+X_j=X_n,n=i+j,i=1,...,2r;j=1,...,s$ in above equation, we have
 	   	  \begin{equation}\label{jy1}
 	   	  	\begin{split}
 	   	  \frac{1}{\lambda^2}g(Ric^{(KerF_*^\perp)}(X)& =  \frac{v}{4}\Big\{(2 r+s+ 2)g_K(X, X) + 3g_K(\omega \mathcal{B}X, X)\Big\}-\frac{1}{4}||3\mathcal{V}[X,X_n]||^2\\&-\frac{\lambda^2}{2}\Big\{2g_K(X, X_n)g_K(\nabla_Xgrad\frac{1}{\lambda^2}, X_n)-g_K(X_n, X_n)\\&g_K\big(\nabla_Xgrad(\frac{1}{\lambda^2}), X\big)-g_K(X, X)g_K(\nabla_{X_n}grad\frac{1}{\lambda^2}, X_n)\Big\}\\&-\frac{\lambda^4}{4}\Big\{\big(g_K(X, X)g_K(X_n, X_n)-(g_K(X, X_n)\big)^2\big)||grad\frac{1}{\lambda^2}||^2\\&+||\big(X(\frac{1}{\lambda^2})X_n - X_n(\frac{1}{\lambda^2})X||^2\big)\Big\},\\
 	\frac{1}{\lambda^2}g(Ric^{(KerF_*^\perp)}(X)   	  &=\dfrac{v}{4}\Big\{(2r+s+2)||X||^2+3g_K(\omega\mathcal{B}X,X)\Big\}-\frac{1}{4}||3\mathcal{V}[X,X_n]||^2\\&-\frac{\lambda^2}{2}\{2g_K(X, X_n)g_K(\nabla_Xgrad\frac{1}{\lambda^2}, X_n)-(2r+s)g_K(\nabla_{X}grad\frac{1}{\lambda^2}, X)\\&
 	   	  -||X||^2g_K(\nabla_{X_n}grad\frac{1}{\lambda^2}, X_n)\big\}-\frac{\lambda^4}{4}\Big\{\big((2r+s)||X||^2\\&
 	   	  -(g_K(X,X_n))^2\big)||grad\frac{1}{\lambda^2}||^2+||\big(X(\frac{1}{\lambda^2})X_n - X_n(\frac{1}{\lambda^2})X\big)||^2\Big\},
 	   	 \end{split}
 	   	  \end{equation}
 	   	hence from \eqref{jy1}, we get the result. 
 	  \end{proof}
 	  
 	 	\section{Conformal pointwise slant Riemannian maps \textit{(CPSRM)} from  Riemannian manifolds to K\"{a}hler manifolds}
 	   In this section, we introduce the notion of conformal pointwise slant Riemannian maps from Riemannian manifolds to K\"{a}hler manifolds with an example and discuss the geometry of such maps.
 	   	  	 
 	   	  	 \begin{definition}
 	   	  	Let $G$ be a conformal Riemannian map from a Riemannian manifold $(L,g_L)$ to a K\"{a}hler manifold $(K,\varphi, g_K)$. If for every point $k\in K$, the Wirtinger angle $\Theta(Z)$ between $\varphi G_*(Z)$ and the space $rangeG_*$ is independent of the choice of $G_*Z$, where $G_*Z\in \Gamma(kerG_*)$ is a nonzero vector, then $G$ is said to be a conformal pointwise slant Riemannian map. In this case, the angle $\Theta$ is regarded as a function on $K,$ known as the slant function of conformal poinwise slant Riemannian map \textit{(CPSRM)}. 
 	   	  	 \end{definition}
 	   	  	 Let, $G$ be a \textit{CPSRM} from a Riemannian manifold $(L,g_L)$ to a K\"{a}hler manifold $(K,\varphi ,g_K)$ with slant function $\Theta$, then for $G_*Z\in\Gamma(rangeG_*),$ we have
 	   	  	 \begin{equation}\label{51}
 	   	  	\varphi G_*Z=\rho G_*Z+\varpi G_*Z,
 	   	  	 \end{equation}
 	   	  	 where $\rho G_*Z\in \Gamma(rangeG_*)$ and $\varpi G_*Z\in \Gamma((rangeG_*)^{\perp}).$ Also for any $P\in \Gamma((rangeG_*)^\perp),$ we have
 	   	  	  \begin{equation}\label{52}
 	   	  	 \varphi P=\mathcal{D}P+\mathcal{E}P,
 	   	  	 \end{equation} 
 	   	  	 where $\mathcal{D}P\in\Gamma(rangeG_*)$ and $\mathcal{E}P\in \Gamma((rangeG_*)^{\perp}).$ Assuming $\eta$ as a orthogonal complementary distribution to $\varpi(\Gamma(rangeG_*))$ in $\Gamma((rangeG_*)^{\perp})$, we can write
 	   	  	 $$\Gamma((rangeG_*)^{\perp})= \varpi(\Gamma(rangeG_*))\oplus \eta.$$ 
 	   	  	 
 	   	  	 \begin{example}
 	   	  	  Consider a Riemannian manifold $(K=\mathbb{R}^4,g_K)$ and a pair of  almost complex structures $\{J_1,J_2\}$ on $K$ satisfying $J_1J_2=-J_2J_1$, where 
 	   	  	 $$J_1(u_1,u_2,u_3,u_4)=(u_3,u_4,-u_1,-u_2)$$
 	   	  	 $$J_2(u_1,u_2,u_3,u_4)=(u_2,-u_1,-u_4,u_3)$$. Let   $t:\mathbb{R}^4\rightarrow\mathbb{R}$ be a real-valued function, hence we can define a complex structure 
 	   	  	 $$J_t=(cost)J_1+(sint)J_2$$ on $K$, then $(K,g_K,J_t)$ is an almost complex structure. Again, consider a map $G:(L=\mathbb{R}^6,g_L)\rightarrow (K=\mathbb{R}^4,g_K)$ from a Riemannian manifold $L$ to a K\"{a}hler manifold $K$, defined by
 	   	  	 $$G(x_1,x_2,x_3,x_4,x_5,x_6)=\pi^a(x_3sinha-x_2cosha,0,x_5cosha-x_4sinha,\sqrt{2}cosb),$$ where $a,b$ are constants. By simple computation, we have
 	   	  	 \begin{equation*}
 	   	  	 (kerG_*)=span\big\{\frac{\partial}{\partial x_1},\pi^acosha\dfrac{\partial}{\partial x_2}+\pi^asinha\dfrac{\partial}{\partial x_3},\pi^acosha\dfrac{\partial}{\partial x_4}+\pi^asinha\dfrac{\partial}{\partial x_5},\dfrac{\partial}{\partial x_6}\big\}
 	   	  	 \end{equation*}
 	   	  	 \begin{equation*}
 	   	  	 (kerG_*)^\perp=span\big\{X=-\pi^acosha\dfrac{\partial}{\partial x_2}+\pi^asinha\dfrac{\partial}{\partial x_3},Y=-\pi^acosha\dfrac{\partial}{\partial x_4}+\pi^asinha\dfrac{\partial}{\partial x_5}\}
 	   	  	 \end{equation*} 
 	   	  	 and
 	   	  	 \begin{equation*}
 	   	  	 rangeG_*=span\big\{G_*X=\pi^{2a}(cosh^2a+sinh^2a,0,0,0),G_*Y=\pi^{2a}(0,0,sinh^2a+cosh^2a,0)\},
 	   	  	 \end{equation*}
 	   	  	hence $G$ is a \textit{CPSRM} from a Riemannian manifold $L$ to a K\"{a}hler manifold $K$ with rank$~G=2,$ $\lambda=\pi^{\alpha}(cosh^2a+sinh^2a)^{1/2}$ and slant function $\Theta=t.$
 	   	  	 \end{example}
 	   	  	 \begin{lemma}
 	   	  		Let $G$ be a conformal Riemannian map from a Riemannian manifold $(L,g_L)$ to a K\"{a}hler manifold $(K,\varphi, g_K)$. Then, $G$ is a \textit{CPSRM} if and only if there exists a constant $\delta \in [-1,0],$ such that
 	   	  		$$\rho^2 G_*(Z)=\delta G_*(Z)$$
 	   	  		for $Z\in\Gamma(kerG_*)^\perp$. If $G$ is a conformal pointwise slant Riemannian map, then $\delta=-cos^2\Theta.$
 	   	  	 \end{lemma}
 	   	  	 By using \eqref{51} and Lemma 4.1, we have following lemma.
 	   	  	 \begin{lemma}
 	   	  	 	Let $G$ be a \textit{CPSRM} from a  Riemannian manifold $(L,g_L)$ to a K\"{a}hler manifold $(K,\varphi,g_K)$ with slant function $\Theta.$ Then, we have
 	   	  	  	 \begin{align}
 	   	  	  	 g_K(\rho G_*(Y),\rho G_*(Z))&=\lambda^2 cos^2\Theta g_L( Y, Z) \label{53}\\
 	   	  	  	 g_K(\varpi G_*(Y),\varpi G_*(Z))&=\lambda^2 sin^2\Theta g_L( Y, Z) \label{54}
 	   	  	  	 \end{align}
 	   	  	  	  for any $Y,Z\in \Gamma((kerF_*)^\perp).$
 	   	  	 \end{lemma}
 	   	  	 \begin{theorem}
 	   	  		Let $G$ be a \textit{CPSRM} from a  Riemannian manifold $(L,g_L)$ to a K\"{a}hler manifold $(K,\varphi,g_K)$. Then any of the following two assertions imply the third one
 	   	  		\begin{enumerate}
 	   	  		\item [$(i)$] $rangeG_*$ is integrable.
 	   	  		\item [$(ii)$] $g_K(\nabla^{G\perp}_Z\varpi\rho G_*Y-\nabla^{G\perp}_Y\varpi\rho G_*Z,P)=g_K(\nabla^{G\perp}_Y\varpi G_*Z-\nabla^{G\perp}_Z\varpi G_*Y,\mathcal{E}P),$
 	   	  		\item[$(iii)$]$g_K((\nabla G_*)(Y,{}^*G_*\mathcal{D}P)^\perp,\varpi G_*Z)=g_K((\nabla G_*)(Z,{}^*G_*\mathcal{D}P)^\perp,\varpi G_*Y),$
 	   	  		 
 	   	  		\end{enumerate}
 	   	  	where,	$Y,Z\in\Gamma((kerG_*)^\perp)$, $P\in \Gamma(range G_*)^\perp)$ and ${}^*G_*$ is adjoint map of $G_*$. 
 	   	  	 \end{theorem}
 	   	  	 \begin{proof}
 	   	  Let $Y,Z\in\Gamma((kerG_*)^\perp)$, $P\in \Gamma((range G_*)^\perp)$ and $K$ is a K\"{a}hler manifold, then we have
 	   	  \begin{equation}\label{54a}
 	   	  	   	  g_K([G_*Y,G_*Z],P)=g_K(\overset{K}{\nabla_Y^G}\varphi G_*Z,\varphi P)-g_K(\overset{K}{\nabla_Z^G}\varphi G_*Y,\varphi P).
 	   	  	   	  \end{equation}
 	   	  Using \eqref{51} in \eqref{54a}, we get
 	   	   \begin{equation}\label{55}
 	   	   	   	\begin{split}
 	   	   	   	 g_K([G_*Y,G_*Z],P)&=g_K(\overset{K}{\nabla_Y^G}\rho G_*Z,\varphi P)+g_K(\overset{K}{\nabla_Y^G}\varpi G_*Z,\varphi P)\\&
 	   	   	   	  -g_K(\overset{K}{\nabla_Z^G}\rho G_*Y,\varphi P)-g_K(\overset{K}{\nabla_Z^G}\rho G_*Y,\varpi P),
 	   	   	   	 \end{split}
 	   	   \end{equation}
 	   	   again, using \eqref{51} and Lemma 4.1 in \eqref{55}, we have
 	   	   \begin{equation*}
 	   	   \begin{split}
 	   	   sin^2\Theta g_K([G_*Y,G_*Z],P)&=-g_K(sin2\Theta Y(\Theta)G_*Z,P)+g_K(sin2\Theta Z(\Theta)G_*Y,P)\\&
 	   	   -g_K(\overset{K}{\nabla_Y^G}\varpi\rho G_*Z, P)
 	   	   +g_K(\overset{K}{\nabla_Z^G}\varpi\rho G_*Y, P)\\&
 	   	   +g_K(\overset{K}{\nabla_Y^G}\varpi G_*Z,\mathcal{D} P) -g_K(\overset{K}{\nabla_Z^G}\varpi G_*Y,\mathcal{D} P)\\&+g_K(\overset{K}{\nabla_Y^G}\varpi G_*Z,\mathcal{E} P) -g_K(\overset{K}{\nabla_Z^G}\varpi G_*Y,\mathcal{E} P),
 	   	   \end{split}
 	   	   \end{equation*}
 	   	   applying \eqref{12} in above equation, then simplifying, we get
 	   	    \begin{equation}\label{56}
 	   	    \begin{split}
 	   	   sin^2\Theta g_K([G_*Y,G_*Z],P)&=
 	   	    g_K(-\nabla_Y^{G^\perp}\varpi\rho G_*Z+\nabla_Z^{G^\perp}\varpi\rho G_*Y, P)\\&+g_K(\nabla_Y^{G^\perp}\varpi G_*Z-\nabla_Z^{G^\perp}\varpi G_*Y,\mathcal{E} P),\\&
 	   	   -g_K(\varpi G_*Z,\overset{K}{\nabla_Y^G}\mathcal{D} P) +g_K(\varpi G_*Y,\overset{K}{\nabla_Z^G}\mathcal{D} P).
 	   	    \end{split}
 	   	    \end{equation}
 	   	Let ${}^*G_{*}$ be the adjoint map of $G_*$, using \eqref{11}, \eqref{13} and \eqref{14}, we get
 	   	 \begin{equation}\label{57}
 	   	 	   	    \begin{split}
 	   	 	   	   sin^2\Theta g_K([G_*Y,G_*Z],P)&=
 	   	 	   	    g_K(-\nabla_Y^{G^\perp}\varpi\rho G_*Z+\nabla_Z^{G^\perp}\varpi\rho G_*Y, P)\\&+g_K(\nabla_Y^{G^\perp}\varpi G_*Z-\nabla_Z^{G^\perp}\varpi G_*Y,\mathcal{E} P)\\&
 	   	 	   	   -g_K(\varpi G_*Z,G_*(\nabla_Y^L{}^*G_{*}\mathcal{D} P)+Y(ln\lambda)\mathcal{D}P
 	   	 	   	   +({}^*G_{*}\mathcal{D}P)(ln\lambda)G_*Y\\&
 	   	 	   	   -g_L(Y,{}^*G_{*}\mathcal{D}P)G_*(grad(ln\lambda))+(\nabla G_*)^\perp(Y,{}^*G_{*}\mathcal{D}P))\\& + g_K(\varpi G_*Y,G_*(\nabla_Z^L{}^*G_{*}\mathcal{D} P)+Z(ln\lambda)\mathcal{D}P
 	   	 	   	   +({}^*G_{*}\mathcal{D}P)(ln\lambda)G_*Z\\&
 	   	 	   	   -g_L(Z,{}^*G_{*}\mathcal{D}P)G_*(grad(ln\lambda))+(\nabla G_*)^\perp(Z,{}^*G_{*}\mathcal{D}P)).\\
 	   	 	   	  & =g_K(-\nabla_Y^{G^\perp}\varpi\rho G_*Z+\nabla_Z^{G^\perp}\varpi\rho G_*Y, P)\\&+g_K(\nabla_Y^{G^\perp}\varpi G_*Z-\nabla_Z^{G^\perp}\varpi G_*Y,\mathcal{E} P)\\&
 	   	 	   	  -g_K((\nabla G_*)^\perp(Y,{}^*G_{*}\mathcal{D}P),\varpi G_*Z)\\& + g_K((\nabla G_*)^\perp(Z,{}^*G_{*}\mathcal{D}P),\varpi G_*Y).
 	   	 	   	    \end{split}
 	   	 	   	    \end{equation}  
 Assuming that assertions $(i)$ and $(ii)$ are true. Applying $(i)$ and $(ii)$ in \eqref{57}, we get $(iii)$. If assertions $(ii)$ and $(iii)$ are true, then from \eqref{57} we have $(i)$. Further, if assertions $(i)$ and $(ii)$ are true, then using $(i)$ and $(ii)$ in \eqref{57}, we get $(ii)$.   	   
 	   	  	 \end{proof}
 	   	  	 
 	   	  	  \begin{theorem}
 	   	  	  	   	  		Let $G$ be a \textit{CPSRM} from a  Riemannian manifold $(L,g_L)$ to a K\"{a}hler manifold $(K,\varphi,g_K)$. Then any two of the following assertions imply the third one
 	   	  	  	   	  		\begin{enumerate}
 	   	  	  	   	  		\item $(rangeG_*)^\perp$ is integrable
 	   	  	  	   	  		\item $(rangeG_*)^\perp$ defines a totally geodesic foliation.
 	   	  	  	   	  		\item $g_K(\overset{K}{\nabla_P}\varpi\rho G_*Y,Q)-g_K(\overset{K}{\nabla_P}\varpi G_*Y,\mathcal{E}Q)=g_K(\overset{K}{\nabla_Q}\varpi\rho G_*Y,P)-g_K(\overset{K}{\nabla_Q}\varpi G_*Y,\mathcal{E}P),$
 	   	  	  	   	  		\end{enumerate}
 	   	  	  	   	  		where $P,Q\in \Gamma ((rangeG_*)^\perp)$ and $Y\in\Gamma((kerG_*)^\perp)$. 
 	   	  	 \end{theorem}
 	   	  	 
 	   	  	 \begin{proof}
 	   	  Let $P,Q\in \Gamma ((rangeG_*)^\perp),$ then for any $Y\in\Gamma((kerG_*)^\perp)$, we have
 	   	  \begin{equation*}
 	   	  g_K([P,Q],G_*Y)=-g_K(\overset{K}{\nabla_P} G_*Y,Q)+g_K(\overset{K}{\nabla_Q} G_*Y,P),
 	   	  \end{equation*}
 	   	  since $K$ is a K\"{a}hler manifold, from above equation, we have
 	   	   \begin{equation} \label{57a}
 	   	   	   	  g_K([P,Q],G_*Y)=-g_K(\overset{K}{\nabla_P} \varphi G_*Y,\varphi Q)+g_K(\overset{K}{\nabla_Q} \varphi G_*Y,\varphi P),
 	   	   	   	  \end{equation}
 	   	   	   	  from \eqref{51}, \eqref{57a} and Lemma 4.1, we have
 	   	   	   	\begin{equation*}
 	   	   	    \begin{split}
 	   	   	   	g_K([P,Q],G_*Y)&=-g_K(\overset{K}{\nabla_P}(cos^2\Theta)G_*Y,Q)+g_K(\overset{K}{\nabla_Q}(cos^2\Theta)G_*Y,P)+g_K(\overset{K}{\nabla_P}\varpi\rho G_*Y,Q)\\&-g_K(\overset{K}{\nabla_Q}\varpi\rho G_*Y,P)-g_K(\overset{K}{\nabla_P}\varpi G_*Y,\varphi Q)+g_K(\overset{K}{\nabla_Q}\varpi G_*Y,\varphi P),
 	   	   	   	\end{split}
 	   	   	   	  \end{equation*}
 	   	   	   using \eqref{52} in above equation, we get
 	   	   	 \begin{equation}\label{57c}
 	   	   	  \begin{split}
 	   	   	  sin^\Theta g_K([P,Q],G_*Y)&=g_K(sin2\Theta P(\Theta)G_*Y,Q)-g_K(sin2\Theta Q(\Theta)G_*Y,P)\\&+g_K(\overset{K}{\nabla_P}\varpi\rho G_*Y,Q)-g_K(\overset{K}{\nabla_Q}\varpi\rho G_*Y,P)\\&-g_K(\overset{K}{\nabla_P}\varpi G_*Y,\mathcal{D} Q+\mathcal{E}Q)+g_K(\overset{K}{\nabla_Q}\varpi G_*Y,\mathcal{D} P+\mathcal{E}P)
 	   	   	  \end{split}
 	   	   	  \end{equation}
 	   	   	  Assuming if the assertions $(i)$ and $(ii)$ are true, from \eqref{57c}, we have $(iii)$. Similarly, assuming assertions $(ii)$ and $(iii)$ , then from \eqref{57c}, we get $(i)$ and further, if assertions $(i)$ and $(iii)$ are true, from \eqref{57c}, we get $(ii)$. Hence, we get the theorem.	  
 	   	  	 \end{proof}
 	   	  	 
 	   	  	 \begin{theorem}
 	   	  	Let $G$ be a \textit{CPSRM} from a  Riemannian manifold $(L,g_L)$ to a K\"{a}hler manifold $(K,\varphi,g_K)$ and $\Theta$ be a proper slant function on $K$. Then  $rangeG_*$ defines a totally geodesic foliation on $K$ if and only if
 	   	  	\begin{equation}
 	   	 g_K\big((\nabla G_*)^\perp(Y,{}^*G_*\mathcal{D}P),\varpi G_*Z\big) =g_K(\nabla_Y^{G\perp}\varpi G_*Z,\mathcal{E}P)-g_K(\nabla_Y^{G\perp}\varpi\rho G_*Z, P),
 	   	  	\end{equation}
 	   	  	where $Y,Z\in\Gamma((kerG_*)^\perp),P\in \Gamma ((rangeG_*)^\perp)$ and ${}^*G_*$ is the adjoint map of $G_*$. 
 	   	  	 \end{theorem}
 	   	  	 \begin{proof}
 	   	  	 Let $Y,Z\in\Gamma((kerG_*)^\perp)$ and $K $ be a K\"{a}hler manifold, then for any $P\in \Gamma ((rangeG_*)^\perp)$, we have
 	   	  	 \begin{equation*}
 	   	  	 g_K(\overset{K}{\nabla_Y^G}G_*Z,P)=g_K(\overset{K}{\nabla_Y^G}\varphi G_*Z,\varphi P),
 	   	  	 \end{equation*}
 	   	  	 applying \eqref{51} in above equation, we get
 	   	  	  \begin{equation}\label{58}
 	   	  	  	   g_K(\overset{K}{\nabla_Y^G}G_*Z,P)=g_K(\overset{K}{\nabla_Y^G}\rho G_*Z+\overset{K}{\nabla_Y^G}\varpi G_*Z,\varphi P),
 	   	  	  \end{equation}
 	   	  	  again, using \eqref{51} in \eqref{58}, we have
 	   	  	  \begin{equation}\label{59}
 	   	  	   	 g_K(\overset{K}{\nabla_Y^G}G_*Z,P)=-g_K(\overset{K}{\nabla_Y^G}\rho^2 G_*Z+\overset{K}{\nabla_Y^G}\varpi\rho G_*Z, P)+g_K(\overset{K}{\nabla_Y^G}\varpi G_*Z,\varphi P).
 	   	  	   \end{equation}
 	   	  	   Using lemma 4.1 and \eqref{12} in \eqref{59}, we can write
 	   	  	   \begin{equation*}
 	   	  	    \begin{split}
 	   	  	    sin^2\Theta g_K(\overset{K}{\nabla_Y^G}G_*Z,P)&=-2sin2\Theta Y(\Theta)g_K(G_*Z,P)-g_K(\nabla_Y^{G\perp}\varpi\rho G_*Z, P)\\&
 	   	  	    -g_K(S_{\varpi G_*Z}G_*Y,\mathcal{D}P)+g_K(\nabla_Y^{G\perp}\varpi G_*Z,\mathcal{E}P),
 	   	  	    \end{split}
 	   	  	    \end{equation*}
 	   	  	    using \eqref{12a} in above equation, further applying \eqref{14}, we get 
 	   	  	     \begin{equation}\label{60}
 	   	  	     \begin{split}
 	   	  	     sin^2\Theta g_K(\overset{K}{\nabla_Y^G}G_*Z,P)&=-g_K(\nabla_Y^{G\perp}\varpi\rho G_*Z, P)+g_K(\nabla_Y^{G\perp}\varpi G_*Z,\mathcal{E}P)\\&
 	   	  	      -g_K\big((\nabla G_*)^\perp(Y,{}^*G_*\mathcal{D}P),\varpi G_*Z\big).
 	   	  	      \end{split}
 	   	  	      \end{equation}
 	   	  	      If $rangeG_*$ defines a totally geodesic foliation on $K$, from \eqref{60} we get the required result.  
 	   	  	 \end{proof}
 	   
 	   	   	  	 \begin{theorem}
 	   	   	  	Let $G$ be a \textit{CPSRM} from a  Riemannian manifold $(L,g_L)$ to a K\"{a}hler manifold $(K,\varphi,g_K)$ and $\Theta$ be a proper slant function on $K$. Then, any two of the following assertions imply the third one
 	   	   	  	\begin{enumerate}
 	   	   	  	\item [$$(i)$$] $(rangeG_*)^\perp$ defines a totally geodesic foliation on $K$,
 	   	   	  	\item [$(ii)$] $G$ is a horizontally homothetic map,
 	   	   	  	\item[$(iii)$] \begin{equation}
 	   	   	  \begin{split}
 	   	   	  	g_K\big((\nabla G_*)^\perp({}^*G_*\mathcal{D}Q,{}^*G_*\mathcal{D}P),\mathcal{E}\varpi G_*Y\big)&=sin^2\Theta\big(\lambda^2g_L(\overset{L}{\nabla}_{{}^*G_*\mathcal{D}Q}{}^*G_*\mathcal{D}P,Y)\\&-g_K(S_{\mathcal{E}P}\mathcal{D}Q,G_*Y)\big)+g_K(\overset{K}{\nabla_P}\varpi G_*Y,\mathcal{E} Q)\\&
 	   	   	 - g_K(\overset{K}{\nabla_P}\varpi\rho G_*Y,Q)-	g_K([P,\mathcal{D}Q],\varpi G_*Y)
 	   	   	  \end{split}
 	   	   	  	\end{equation}
 	   	   	  	\end{enumerate}   
 	   	   	  	where $P,Q\in \Gamma ((rangeG_*)^\perp),Y\in\Gamma((kerG_*)^\perp)$ and ${}^*G_*$ is the adjoint map of $G_*$. 
 	   	   	  	 \end{theorem}
 	   	   	  	 \begin{proof}
 	   	   	  	 Let $P,Q\in \Gamma ((rangeG_*)^\perp),Y\in\Gamma((kerG_*)^\perp)$ and $K$ be a K\"{a}hler manifold, from \eqref{51} we can write
 	   	   	  	 \begin{equation}\label{61}
 	   	   	  	 g_K(\overset{K}{\nabla_P}Q,G_*Y)=-g_K(\overset{K}{\nabla_P}\rho G_*Y,\varphi Q)-g_K(\overset{K}{\nabla_P}\varpi G_*Y,\varphi Q),
 	   	   	  	 \end{equation}
 	   	   	  	 using \eqref{51} and \eqref{52} in \eqref{61}, we get
 	   	   	  	 \begin{equation*}
 	   	   	  	  g_K(\overset{K}{\nabla_P}Q,G_*Y)=g_K(\overset{K}{\nabla_P}\rho^2 G_*Y, Q)+g_K(\overset{K}{\nabla_P}\varpi\rho G_*Y, Q)+g_K(\overset{K}{\nabla_P}\mathcal{D}Q,\varpi G_*Y)-g_K(\overset{K}{\nabla_P}\varpi G_*Y,\mathcal{E} Q),
 	   	   	  	 \end{equation*}
 	   	   	  	 applying lemma 4.1 in above equation and simplifying, we get
 	   	   	  	  \begin{equation}\label{62}
 	   	   	  	 \begin{split}
 	   	   	  	  g_K(\overset{K}{\nabla_P}Q,G_*Y)&=cos^2\Theta g_K(\overset{K}{\nabla_P}Q,G_*Y)+g_K(\overset{K}{\nabla_P}\varpi\rho G_*Y, Q)
 	   	   	  	  +g_K(\overset{K}{\nabla_{\mathcal{D}Q}}\mathcal{E}P,\varphi\varpi G_*Y)\\&
 	   	   	  	  +g_K(\overset{K}{\nabla_{\mathcal{D}Q}}\mathcal{D}P,\varphi\varpi G_*Y)+g_K([P,\mathcal{D}Q],\varpi G_*Y)-g_K(\overset{K}{\nabla_P}\varpi G_*Y,\mathcal{E} Q).
 	   	   	  	 \end{split}
 	   	   	  	  \end{equation}
 	   	   	  	 Suppose ${}^*G_*$ be the adjoint map of $G_*$, using \eqref{11}, \eqref{12}, \eqref{13} and \eqref{14} in \eqref{62}, we get
 	   	   	  	   \begin{equation*}
 	   	   	  	   	  \begin{split}
 	   	   	  	   	   sin^2\Theta g_K(\overset{K}{\nabla_P}Q,G_*Y)&=g_K(\overset{K}{\nabla_P}\varpi\rho G_*Y, Q)-g_K(S_{\mathcal{E}P}\mathcal{D}Q,\varphi \varpi G_*Y)
 	   	   	  	   	   +g_K(\nabla_{\mathcal{D}Q}^{G\perp}\mathcal{E}P,\varphi\varpi G_*Y)\\&
 	   	   	  	   	    +g_K(\overset{L}{\nabla}_{{}^*G_*\mathcal{D}Q}{}^*G_*\mathcal{D}P+{}^*G_*\mathcal{D}Q(ln\lambda)\mathcal{D}P+{}^*G_*\mathcal{D}P(ln\lambda)\mathcal{D}Q\\&-g_L({}^*G_*\mathcal{D}Q,{}^*G_*\mathcal{D}PG_*(grad(ln\lambda))+(\nabla G_*)^\perp({}^*G_*\mathcal{D}Q,{}^*G_*\mathcal{D}P),\varphi\varpi G_*Y)\\&
 	   	   	  	   	    +g_K([P,\mathcal{D}Q],\varpi G_*Y)-g_K(\overset{K}{\nabla_P}\varpi G_*Y,\mathcal{E} Q),
 	   	   	  	   	 \end{split}
 	   	   	  	   	  \end{equation*}
 	   	   	  	   	 applying \eqref{52} and lemma 4.2 in above equation and simplifying, we get
 	   	   	  	   	  \begin{equation}\label{63}
 	   	   	  	   	  \begin{split}
 	   	   	  	   	  sin^2\Theta g_K(\overset{K}{\nabla_P}Q,G_*Y)   &=g_K(\overset{K}{\nabla_P}\varpi\rho G_*Y, Q)+(\nabla G_*)^\perp({}^*G_*\mathcal{D}Q,{}^*G_*\mathcal{D}P),\mathcal{E}\varpi G_*Y)\\& +g_K([P,\mathcal{D}Q],\varpi G_*Y)-g_K(\overset{K}{\nabla_P}\varpi G_*Y,\mathcal{E} Q)+ sin^2\Theta \Big (g_K(S_{\mathcal{E}P}\mathcal{D}Q,G_*Y)\\&
 	   	   	  	   	  -\lambda^2g_L(\mathcal{H}\overset{L}{\nabla}_{{}^*G_*\mathcal{D}Q}{}^*G_*\mathcal{D}P,Y)-{}^*G_*\mathcal{D}Q(ln\lambda)g_K(\mathcal{D}P,G_*Y)\\&-{}^*G_*\mathcal{D}P(ln\lambda)g_K(\mathcal{D}Q,G_*Y)+Y(ln\lambda)g_K(\mathcal{D}P,\mathcal{D}Q)\Big ).
 	   	   	  	   	  \end{split}
 	   	   	  	   	  \end{equation}
 	   	   	  	   	  Suppose, assertions $(i)$ and $(ii)$ are true, then from \eqref{63} we have $(iii)$. Similarly, if $(ii)$ and $(iii)$ are true, from \eqref{63} we get assertion $(i)$. Again, if $(i)$ and $(iii)$ are true, from \eqref{63}, we have
 	   	   	  	   	  $$sin^2\Theta\big(-{}^*G_*\mathcal{D}Q(ln\lambda)g_K(\mathcal{D}P,G_*Y)-{}^*G_*\mathcal{D}P(ln\lambda)g_K(\mathcal{D}Q,G_*Y)+Y(ln\lambda)g_K(\mathcal{D}P,\mathcal{D}Q)\big )=0,$$ 
 	   	   	  	   	since $\Theta\neq0$ and putting $P=Q$, we have
 	   	   	  	   	$$2g_K\big(\mathcal{D}P,G_*(grad(ln\lambda))\big)g_K(\mathcal{D}P,G_*Y)=g_K\big(G_*X,G_*(grad(ln\lambda))\big)g_K(\mathcal{D}P,\mathcal{D}P),$$ 
 	   	   	  	   	which is possible only if $\lambda$ is constant on $(kerG_*)^\perp$, hence the assertion $(ii)$.  
 	   	   	  	 \end{proof} 
 	   	   	  	 \begin{theorem}
 	   	   	  		Let $G$ be a \textit{CPSRM} from a  Riemannian manifold $(L,g_L)$ to a K\"{a}hler manifold $(K,\varphi,g_K)$. Then, $G$ is harmonic if the following conditions are satisfied:
 	   	   	  		\begin{enumerate}
 	   	   	  		\item [$(i)$] \begin{equation}
 	   	   	  		\begin{split}
 	   	   	  		&trace\big\{S_{\varpi\rho G_*(\cdot)}G_*(\cdot)+S_{\mathcal{E}\varpi G_*(\cdot)}G_*(\cdot)+G_*\big(\nabla^L_{(\cdot)}{}^*G_*(\mathcal{D}\varpi G_*(\cdot))\big)\\&-(\nabla G_*)\big(\cdot,{}^*G_*(\mathcal{D}\varpi G_*(\cdot))\big)^{rangeG_*}-sin2\Theta (\cdot)(\Theta)G_*(\cdot)-sin^2\Theta G_*(\nabla^L_{(\cdot)}(\cdot))\big\}=0,
 	   	   	  		\end{split}
 	   	   	  		\end{equation}
 	   	   	  		
 	   	   	  		\item[$(ii)$] \begin{equation}
 	   	   	  		trace\big\{\nabla_{(\cdot)}^{G\perp}\varpi\rho G_*(\cdot)+(\nabla G_*)^\perp\big(\cdot,{}^*G_*(\mathcal{D}\varpi G_*(\cdot))\big)+\nabla_{(\cdot)}^{G\perp}\mathcal{E}\varpi G_*(\cdot)\big\}=0,
 	   	   	  		\end{equation}
 	   	   	  		\item[$(iii)$] fibers are minimal.
 	   	   	  		\end{enumerate}
 	   	   	  	 \end{theorem}
 	   	   	  	 \begin{proof}
 	   	   	  	 Let $Y\in\Gamma(kerF_*)^\perp$ and $K$ be a K\"{a}hler manifold, from \eqref{11} we have
 	   	   	  	 \begin{equation}\label{71}
 	   	   	  	 (\nabla G_*)(Y,Y)=-\varphi\overset{K}{\nabla_Y^G}\varphi G_*Y-G_*(\nabla_Y^MY),
 	   	   	  	 \end{equation}
 	   	   	  	 further, from \eqref{51}, \eqref{52} and \eqref{71}, we get
 	   	   	  	 \begin{equation*}
 	   	   	  	(\nabla G_*)(Y,Y)=-\overset{K}{\nabla_Y^G}\rho^2 G_*Y-\overset{K}{\nabla_Y^G}\varpi\rho G_*Y-\overset{K}{\nabla_Y^G}\mathcal{D}\varpi G_*Y-\overset{K}{\nabla_Y^G}\mathcal{E}\varpi G_*Y-G_*(\overset{L}{\nabla_Y}Y),
 	   	   	  	 \end{equation*}
 	   	   	  	 using lemma 4.1 and \eqref{12a} in above equation, we have
 	   	   	  	\begin{equation}\label{72}
 	   	   	  	\begin{split}
 	   	   	  (\nabla G_*)(Y,Y)&=cos^2\Theta\overset{K}{\nabla_Y^G} G_*Y-sin2\Theta Y(\Theta)G_*Y+S_{\varpi\rho G_*Y}G_*Y -\nabla_Y^{G\perp}\varpi\rho G_*Y\\&-\overset{K}{\nabla_Y^G}\mathcal{D}\varpi G_*Y+S_{\mathcal{E}\varpi G_*Y}-\nabla_Y^{G\perp}\mathcal{E}\varpi G_*Y-G_*(\overset{L}{\nabla_Y}Y),
 	   	   	  	\end{split}
 	   	   	  	\end{equation}
 	   	   	  	Let ${}^*G_*$ be the adjoint map of $G_*$, using \eqref{11} in \eqref{72} and separating components of $rangeG_*$ and $(rangeG_*)^\perp$, we get
 	   	   	  	\begin{equation}\label{73}
 	   	   	  	\begin{split}  	   	  		
 	   	   	  	 sin^2\Theta(\nabla G_*)(Y,Y)^{rangeG_*}&= S_{\varpi\rho G_*Y}G_*Y+S_{\mathcal{E}\varpi G_*Y}G_*Y+G_*\big(\overset{L}{\nabla_{Y}}{}^*G_*(\mathcal{D}\varpi G_*Y)\big)\\&-(\nabla G_*)\big(Y,{}^*G_*(\mathcal{D}\varpi G_*Y)\big)^{rangeG_*}-sin2\Theta Y(\Theta)G_*Y\\&
 	   	   	  	-sin^2\Theta G_*(\overset{}{\nabla_{Y}}Y),
 	   	   	  	\end{split}
 	   	   	  	\end{equation} 
 	   	   	and 
 	   	   	\begin{equation}\label{74}
 	   	  	sin^2\Theta(\nabla G_*)^{\perp}(Y,Y)=-\nabla_{Y}^{G\perp}\varpi\rho G_*Y-(\nabla G_*)^\perp\big(Y,{}^*G_*(\mathcal{D}\varpi G_*Y)\big)-\nabla_{Y}^{G\perp}\mathcal{E}\varpi G_*Y,
 	   	  	\end{equation}
 	   	  	Again, let $W\in\Gamma(kerG_*)$, from \eqref{11}, we have
 	   	  	\begin{equation}\label{75}
 	   	  (\nabla G_*)(W,W)=-G_*(\nabla_WW)=-G_*(\mathcal{T}_WW)
 	   	  	\end{equation} 
 	   	  	Thus, from \eqref{73}, \eqref{74} and \eqref{75}, we get the required results. 	 	   	   	  		
 	   	   	  	 \end{proof}
 	   	   	  	 \begin{theorem}
 	   	   	  		Let $G$ be a \textit{CPSRM} from a  Riemannian manifold $(L,g_L)$ to a K\"{a}hler manifold $(K,\varphi,g_K)$ and $\Theta$ be the slant function on $K$. Then we have
 	   	   	  		\begin{equation}
 	   	   	  	\begin{split}
 	   	   	  		sin^4\Theta||(\nabla G_*)(Y,Y)^{rangeG_*}||^2&\geq sin^4\Theta ||G_*(\overset{L}{\nabla_Y}Y)||^2+||S_{\mathcal{E}\varpi G_*Y}G_*Y||^2+||S_{\varpi\rho G_*Y}G_*Y||^2\\&
 	   	   	  		+||G_*\big(\overset{L}{\nabla_{Y}}{}^*G_*(\mathcal{D}\varpi G_*Y)\big)||^2+||(\nabla G_*)\big(Y,{}^*G_*(\mathcal{D}\varpi G_*Y)\big)^{rangeG_*}||^2\\&
 	   	   	  		+2\Big\{sin^2\Theta\big(\lambda^2sin2\Theta Y(\Theta)g_L(\nabla_YY,Y)-g_L(\overset{L}{\nabla_Y}Y,\overset{L}{\nabla_{Y}}{}^*G_*(\mathcal{D}\varpi G_*Y))\\&
 	   	   	  		-g_K\big(S_{\varpi\rho G_*Y}G_*Y
 	   	   	  		+S_{\mathcal{E}\rho G_*Y}G_*Y\\&
 	   	   	  		-(\nabla G_*)(Y,{}^*G_*(\mathcal{D}\varpi G_*Y))^{rangeG_*},G_*(\overset{L}{\nabla_Y}Y)\big)\big)\\&
 	   	   	  		-sin2\Theta Y(\Theta)\big(\lambda^2g_L(\overset{L}{\nabla_Y}{}^*G_*(\mathcal{D}\varpi G_*Y),Y)-g_K(S_{\varpi\rho G_*Y}G_*Y)\\&
 	   	   	  	    +S_{\mathcal{E}\varpi G_*Y}G_*Y-(\nabla G_*)\big(Y,{}^*G_*(\mathcal{D}\varpi G_*Y)\big)^{rangeG_*},G_*Y\big)\\&
 	   	   	  		+g_K\big(S_{\mathcal{E}\varpi G_*Y}G_*Y+G_*(\overset{L}{\nabla_Y}{}^*G_*(\mathcal{D}\varpi G_*Y))\\&
 	   	   	    	-(\nabla G_*)\big(Y,{}^*G_*(\mathcal{D}\varpi G_*Y)\big)^{rangeG_*},S_{\varpi\rho G_*Y}G_*Y\big)+g_K\big(S_{\mathcal{E}\varpi G_*Y}G_*Y\\&
 	   	   	    	-(\nabla G_*)\big(Y,{}^*G_*(\mathcal{D}\varpi G_*Y)\big)^{rangeG_*},G_*(\overset{L}{\nabla_Y}{}^*G_*(\mathcal{D}\varpi G_*Y))\big)\Big\},
 	   	   	  		\end{split}
 	   	   	  		\end{equation}
 	   	   equality holds when $\Theta$ is constant, also
 	   	   \begin{equation}
 	   	   \begin{split}
 	   	   sin^4\Theta||(\nabla G_*)^\perp(Y,Y)||^2&=||\nabla_Y^{G\perp}\varpi\rho G_*Y||^2 +||\nabla_Y^{G\perp}\mathcal{E}\varpi G_*Y||^2 +||(\nabla G_*)^\perp(Y,{}^*G_*(\mathcal{D}\varpi G_*Y))||^2\\& +2\Big\{g_K\big(\nabla_Y^{G\perp}\mathcal{E}\varpi G_*Y+(\nabla G_*)^\perp(Y,{}^*G_*(\mathcal{D}\varpi G_*Y)),\nabla_Y^{G\perp}\varpi\rho G_*\big)\\&
 	   	   +g_K\big(\nabla_Y^{G\perp}\mathcal{E}\varpi G_*Y,\nabla_Y^{G\perp}\mathcal{E}\varpi G_*Y\big)\Big\}
 	   	   \end{split}
 	   	   \end{equation}	  		
 	   	   	  	 \end{theorem}
 	   	   	  	 \begin{proof}
 	   	   	  	After taking the products of \eqref{73} and \eqref{74} by itself and further rearranging terms we get the required results. 
 	   	   	  	 \end{proof}
 	   	   	  	 
 	   	   	  	\section{Acknowledgments}
 	   	   	   The first author is thankful to UGC for providing financial assistance in terms of MANF scholarship vide letter with UGC-Ref. No. 1844/(CSIR-UGC NET JUNE 2019). The second author is thankful to DST Gov. of India for providing financial support in terms of DST-FST label-I grant vide sanction number SR/FST/MS-I/2021/104(C).

 	\end{document}